\documentclass{article}
\usepackage{amsmath,amssymb,amsthm,dsfont,hyperref}
\begin{document}
    \title{On VOAs Associated to Jordan Algebras of Type $C$}
    \author{Hongbo Zhao}
    \date{}
    \maketitle
    \allowdisplaybreaks
    \newtheorem{theorem}{Theorem}
    \newtheorem{definition}{Definition}
    \newtheorem{example}{Example}
    \newtheorem{lemma}{Lemma}
    \newtheorem{proposition}{Proposition}
    \section{Introduction and Main Results}
    \par{
        In is well known that for any $\mathbb{Z}_{\geq 0}$-graded vertex operator algebra (VOA) such that $V_0=\mathbb{C}1,V_1=\{0\}$, $V_2$ has a commutative (but not necessarily associative) algebra structure, with the operation given by $a\circ b=a(1)b$. This algebra $V_2$ is called the Griess algebra of $V$.
    }
    \par{
         The case when the Greiss algebra $V_2$ is isomorphic to a finite dimensional simple Jordan algebra is of particular interest to us.  In \cite{Lam1} and \cite{Lam2}, Lam constructed VOAs whose Griess algebras are simple Jordan algebras of Hermitian type. In \cite{AM}, Ashihara and Miyamoto constructed a family of VOAs $V_{\mathcal{J},r}$ parametrized by $r\in\mathbb{C}$, whose Griess algebras are isomorphic to type $B$ simple Jordan algebras $\mathcal{J}$.
    }
    \par{
        It is well known in the theory of VOA that we have many `universal constructions'. The most important examples include the level $k$ universal VOA $V^k(\mathfrak{g})$ associated to a finite dimensional simple Lie algebra $\mathfrak{g}$, and the universal Virasoro VOA $M(k,0)$ with central charge equals $k$. It is also important to study the irreducibility of these VOAs, and the cases when $V^k(\mathfrak{g})$ and $M(k,0)$ are reducible are of particular interest, which is sometimes called `degenerate'. Let $V_k(\mathfrak{g})$ and $L(k,0)$ be the VOA simple quotients of $V^k(\mathfrak{g})$ and $M(k,0)$ respectively, then in the `degenerate' cases, these simple VOAs have many interesting properties. The irreducibility of the universal VOAs $V^k(\mathfrak{g})$ and $M(k,0)$, together with the constructions and properties of the simple VOAs $V_k(\mathfrak{g})$ and $L(c,0)$, has already been extensively studied in literatures for decades.
    }
    \par{
        In our previous paper \cite{Z2}, we made a further study of the VOA $V_{\mathcal{J},r}$ constructed by Ashihara and Miyamoto. We showed that the VOAs $V_{\mathcal{J},r}$ are the `universal VOAs' satisfying $V_0=\mathbb{C}1,V_1=\{0\}$, $V_2\simeq \mathcal{J}$, which resembles the VOAs $V^k(\mathfrak{g})$ and $M(c,0)$ mentioned above. We also explicitly constructed the simple quotients $\bar{V}_{\mathcal{J},r}$ for $r\in\mathbb{Z}_{\neq 0}$, using dual-pair type constructions. We also reprove that $V_{\mathcal{J},r}$ is simple (or equivalently, $V_{\mathcal{J},r}=\bar{V}_{\mathcal{J},r}$) if and only if $r\notin \mathbb{Z}$ using a different method. It follows that the VOA constructed by Lam for a type $B$ simple Jordan algebra, is actually isomorphic to $\bar{V}_{\mathcal{J},1}$.
    }
    \par{
        The main result of this paper, is to extend the construction of Ashihara and Miyamoto for the VOA $V_{\mathcal{J},r}$, to an arbitrary Hermitian type Jordan algebra $\mathcal{J}$. We study the simplicities of these VOAs $V_{\mathcal{J},r}$, and prove similar irreducibility results. We construct the corresponding simple quotients $\bar{V}_{\mathcal{J},r}$ explicitly using dual-pair type constructions, which are parallel to the results in \cite{Z2}. We also show that any simple VOA constructed by Lam in \cite{Lam2}, is isomorphic to one of the simple quotients $\bar{V}_{\mathcal{J},1}$.
    }
    \par{
        In particular, we have
        \begin{theorem}
            \begin{enumerate}
                Let $\mathcal{J}=\mathcal{J}(W)$ be the type $C$ Jordan algebra associated to a symplectic space $W$, $dim(W)\geq 4$, then
                \item For arbitrary complex number $r$, there is a VOA $V=V_{\mathcal{J},r}$ such that $V_0=\mathbb{C}1$, $V_1=\{0\}$. For the Greiss algebra $V_2$ we have
                $$
                    V_2\simeq \mathcal{J},
                $$
                and the central charge of $V_{\mathcal{J},r}$ equals $-dim(W)r$.
                \item $V=V_{\mathcal{J},r}$ is generated by $V_2$.
                \item $V_{\mathcal{J},r}$ is simple if and only if $r\notin\mathbb{Z}$.
            \end{enumerate}
        \end{theorem}
    }
    \par{
        The content of this article is divided into four parts. In Section 2 we briefly review some basic facts about finite dimensional simple Jordan algebras. In Section 3 we review the VOA $V_{\mathcal{J},r}$, where $\mathcal{J}$ is a Type $B$ Jordan Algebra, and our exposition is a little bit different from the one given in \cite{AM}. It will be seen that neither the central extension, or normal ordering is needed in the construction. Similar to the approach given in Section 3, we construct the VOA $V_{\mathcal{J},r}$, where $\mathcal{J}$ is a Type $C$ Jordan Algebra in Section 4, and we prove (1) and (2) in Theorem 1, which is similar to the main results in \cite{AM} and \cite{NS}. Finally in Section 5, we prove the simplicity result for the VOA $V_{\mathcal{J},r}$, where $\mathcal{J}$ is a Type $C$ Jordan Algebra.
    }
    \section{Finite Dimensional Simple Jordan Algebras and Jordan Frames}
    \par{
        In this section, we recall some basic facts about finite dimensional simple Jordan algebras.
    }
    \par{
        Recall that a Jordan algebra $\mathcal{J}$ is a vector space, together with a bilinear map $\circ:\mathcal{J}\times \mathcal{J}\rightarrow \mathcal{J}$ called Jordan product, satisfying
            \begin{align*}
               &x\circ y=y\circ x,\\
               &(x\circ y)\circ(x\circ x)=x\circ (y\circ (x\circ x)),
            \end{align*}
        for all $x,y\in\mathcal{J}$. It follows from the definition that a Jordan algebra is commutative, but it is not necessarily associative.
    }
    \par{
        For any associative algebra $A$, it has a structure of Jordan algebra, with the Jordan product given by:
            $$
                x\circ y=\frac{1}{2}(xy+yx),\text{ for all } x,y\in A.
            $$
            We call $\{\mathcal{J},\circ\}$ a special Jordan algebra if it is isomorphic to a Jordan subalgebra of $\{A,\circ\}$, for some associative algebra $A$. Otherwise we say the Jordan algebra $\{\mathcal{J},\circ\}$ is exceptional.
    }
    \par{
        It is well known that the finite dimensional simple Jordan algebras over an algebraically closed field of characteristic zero is classified by A. Albert in \cite{A}. When the ground field is $\mathbb{C}$, there are five types of simple Jordan algebras, called type $A$, $B$, $C$, $D$ and $E$ Jordan algebras respectively. For details about the general theory of Jordan algebras, the readers can consult \cite{Mc}, \cite{FK94}.
    }
    \par{
        For our purpose, we only describe type $A$, $B$ and $C$ Jordan algebras here. We realize type $A$ and type $B$ Jordan algebras using tensors. Let $(\mathfrak{h},(\cdot,\cdot))$ be a $d$-dimensional vector space with a non-degenerate bilinear form $(\cdot,\cdot)$. Then $\mathfrak{h}\otimes \mathfrak{h}$ has an associative algebra structure:
        $$
            (a\otimes b)(u\otimes v)=(b,u)a\otimes v,
        $$
        which induces a Jordan algebra structure on $\mathfrak{h}\otimes \mathfrak{h}$:
        \begin{align}
            x\circ y=\frac{1}{2}(xy+yx),\text{ for all } x,y\in \mathfrak{h}\otimes \mathfrak{h}.\label{Jordan}
        \end{align}
        We call $\mathfrak{h}\otimes \mathfrak{h}$ the type $A$ Jordan algebra associated to $\mathfrak{h}$, denoted by $\mathcal{J}_A(\mathfrak{h})$, and we set
        \begin{align}
            \tilde{L}_{a,b}:=a\otimes b \in \mathfrak{h}\otimes \mathfrak{h}. \label{JordanA}
        \end{align}
        We further assume that the bilinear form $(\cdot,\cdot)$ is symmetric, and let $\mathcal{J}_B(\mathfrak{h})$ be the Jordan subalgebra of $\mathcal{J}_A(\mathfrak{h})$, which consists of symmetric square tensors:
        \begin{align}
            \mathcal{J}_B(\mathfrak{h}):=S^2(\mathfrak{h})=span\{L_{a,b}|\,a,b\in\mathfrak{h}\},\;L_{a,b}:=a\otimes b+b\otimes a. \label{JordanB}
        \end{align}
        We call $\mathcal{J}_B(\mathfrak{h})$ the type $B$ Jordan algebra associated to $\mathfrak{h}$. It is easy to check that these definitions coincide with the definitions using matrices.
    }
    \par{
        We realize type $C$ Jordan algebras in a similar way. Let $W$ be a $2d$-dimensional space with a non-degenerated bilinear form $\langle \cdot,\cdot \rangle$, then we have the corresponding type $A$ Jordan algebra $\mathcal{J}_A(W)$. We further assume that $\langle\cdot,\cdot\rangle$ is skew-symmetric, and we consider the Jordan subalgebra $\mathcal{J}_C(W)$ of $\mathcal{J}_{A}(W)$, which consists of anti-symmetric tensors:
        \begin{align}
            \mathcal{J}_C(W):=\wedge^2(W)=span\{L_{a,b}|\,a,b\in W\},\;L_{a,b}:=a\otimes b-b\otimes a. \label{JordanC}
        \end{align}
        We call $\mathcal{J}_C(W)$ the type $C$ Jordan algebra associated to $W$. It is also easy to check that this definition is the same as the definition of type $C$ Jordan algebras using skew-symmetric matrices.
    }
    \par{
        By (\ref{Jordan}), (\ref{JordanA}), (\ref{JordanB}), and (\ref{JordanC}), it is easy to check the following explicit formulas about Jordan products in type $A$, $B$ and $C$ Jordan algebras, which are useful in later discussions:
        \begin{align}
            &L_{a,b}\circ L_{u,v}=\frac{1}{2}(b,u)L_{a,v}+\frac{1}{2}(b,v)L_{a,u}+\frac{1}{2}(a,u)L_{b,v}+\frac{1}{2}(a,v)L_{b,u},\notag\\
            &\text{ for all } L_{a,b},L_{u,v}\in \mathcal{J}_B(\mathfrak{h});\notag\\
            &L_{a,b}\circ L_{u,v}=\frac{1}{2}(b,u)L_{a,v}-\frac{1}{2}(b,v)L_{a,u}-\frac{1}{2}(a,u)L_{b,v}+\frac{1}{2}(a,v)L_{b,u},\label{JordanCmul}\\
            &\text{ for all } L_{a,b},L_{u,v}\in \mathcal{J}_C(W);\notag\\
            &\tilde{L}_{a,b}\circ \tilde{L}_{u,v}=\frac{1}{2}(b,u)\tilde{L}_{a,v}+\frac{1}{2}(a,v)\tilde{L}_{b,u},\notag\\
            &\text{ for all } \tilde{L}_{a,b},\tilde{L}_{u,v}\in \mathcal{J}_A(\mathfrak{h}).\notag
        \end{align}
    }
    \par{
        We write down the Jordan frame of type $B$ and type $C$ Jordan algebras for later use. Recall that the Jordan frame of a simple Jordan algebra \cite{FK94} is a set of idempotents $u_1,\cdots u_n$ such that each $u_i$ is not a sum of two non-zero mutually orthogonal idempotents, and
        $$
            \sum u_i=e,\,u_i\circ u_j=\delta_{i,j}u_i
        $$
        for all $i,j$, where $e$ is the identity element in the Jordan algebra. For type $B$ and type $C$ Jordan algebras it is easy to find out the corresponding Jordan frames. Let $\{e_1,\cdots,e_d\}$ be an orthonormal basis of $\mathfrak{h}$, and $\{\psi_1,\cdots ,\psi_n,\psi_1^{*},\cdots ,\psi_n^{*}\}$ be a symplectic basis of $W$ such that
        $$
            \langle \psi^{*}_i,\psi_j\rangle=\delta_{i,j},\,\langle \psi^{*}_i,\psi^{*}_j\rangle=\langle \psi_i,\psi_j\rangle=0
        $$
        for all $1\leq i,j\leq n$. Then we check that the Jordan frame of $\mathcal{J}_B(\mathfrak{h})$ can be given by
        $$
            \{\frac{1}{2}L_{e_i,e_i}|\,i=1,\cdots,d\},
        $$
        and for the Jordan algebra $\mathcal{J}_C(W)$:
        $$
            \{L_{\psi_i,\psi_i^{*}}|\,i=1,\cdots,n\}.
        $$
        It will be seen later that the elements in the Jordan frame for each case corresponds to a set of mutually orthogonal Virasoro elements in the Greiss algebras, which will play a role in analyzing the corresponding VOA structures.
    }
    \section{Construction of the VOA $V_{\mathcal{J},r}$, Where $\mathcal{J}$ is a Type $B$ Jordan Algebra}
    \par{
        In this section we review the VOA $V_{\mathcal{J},r}$ constructed by Ashihara and Miyamoto, Where $\mathcal{J}$ is a Type $B$ Jordan Algebra. We explain the corresponding Lie algebra $\mathcal{L}$ as a central extension of another familiar infinite dimensional Lie algebra. This approach is slightly different from the original one given in \cite{AM}, but they are essentially the same. Our approach is very similar to the one given by Kac and Radul in \cite{KR1}, and it is easy to see the analogy between $V_{\mathcal{J},r}$ and the universal VOA $M_c$ given in \cite{KR1}. The construction described in this section will also inspire construction for type $C$ Jordan algebra, which will be given in Section 3.
    }
    \par{
        We fix the notation that $\mathfrak{h}$ denotes a $d$-dimensional vector space with a non-degenerate symmetric bilinear form $(\cdot,\cdot)$, $\{e_1,\cdots,e_d\}$ be an orthonormal basis of $\mathfrak{h}$, and $\mathcal{J}$ be the type $B$ Jordan algebra $\mathcal{J}_B(\mathfrak{h})$. We recall the following infinite dimensional Lie algebra $\hat{\mathfrak{h}}$ associated to $\mathfrak{h}$:
        \begin{align*}
            \hat{\mathfrak{h}}=\mathfrak{h}\otimes\mathbb{C}[t,t^{-1}]\oplus \mathbb{C}c.
        \end{align*}
        The Lie bracket over $\hat{\mathfrak{h}}$ is given by:
        \begin{align*}
            [a(m),b(n)]=m(a,b)\delta_{m+n,0}c,\;\;
            [x,c]=0,\text{ for all } x\in \hat{\mathfrak{h}}.\notag
        \end{align*}
        It is well known that
        $$
            \hat{\mathfrak{h}}_{-}\stackrel{def.}{=}\mathfrak{h}\otimes\mathbb{C}t^{-1}[t^{-1}]
        $$
        is a commutative Lie subalgebra of $\hat{\mathfrak{h}}$. The Fock space $S(\hat{\mathfrak{h}}_{-})\simeq U(\hat{\mathfrak{h}}_{-})\cdot 1$ is a left $U(\hat{\mathfrak{h}})$-module, and $S(\hat{\mathfrak{h}}_{-})$ has a vertex operator algebra structure \cite{FLM}. We denote this VOA by
        $\mathcal{H}(\mathfrak{h})$.
    }
    \par{
        Given a vector space $W$ with an antisymmetric bilinear form $\langle \cdot,\cdot\rangle$, the symmetric square $S^2(W)$ is a Lie algebra, and for $ab,uv\in S^2(W)$, the Lie bracket is given by
        \begin{align}
            [ab,uv]=\langle b,u \rangle av+\langle a,u \rangle bv+\langle b,v \rangle ua+\langle a,v \rangle ub. \label{bracket}
        \end{align}
        In particular, if $W$ is a finite dimensional symplectic space, the symmetric square $S^2(W)$ is isomorphic to the finite dimensional symplectic Lie algebra $\mathfrak{sp}(W)$. The corresponding Lie bracket is the same as (\ref{bracket}), and this construction is essentially related to the Oscillator representations.
    }
     \par{
        We apply this to the case
        $$
            \widetilde{W_{\infty}}{=}\mathfrak{h}\otimes\mathbb{C}[t,t^{-1}].
        $$
        It is obvious that there is an anti-symmetric bilinear form $\langle\cdot,\cdot\rangle'$ on the space of Laurant polynomials $\mathbb{C}[t,t^{-1}]$
        $$
            \langle t^m,t^n\rangle'=m\delta_{m+n,0},
        $$
        therefore $\widetilde{W_{\infty}}$ also has an anti-symmetric bilinear form $\langle\cdot,\cdot\rangle=(\cdot,\cdot)\otimes \langle\cdot,\cdot\rangle'$
        $$
            \langle a(m),b(n)\rangle=m(a,b)\delta_{m+n,0},
        $$
        which is the same as the value of $[a(m),b(n)]$ by taking $c=1$. Hence $S^2(\widetilde{W_{\infty}})$ is also a Lie algebra with the bracket given by (\ref{bracket}). We note that $\langle\cdot,\cdot\rangle$ restricted to the following subspace
        $$
            W_{\infty}\stackrel{def.}{=}\mathfrak{h}\otimes\mathbb{C}t[t]\bigoplus \mathfrak{h}\otimes\mathbb{C}t^{-1}[t^{-1}]
        $$
        is non-degenerate, and the Lie subalgebra $S^2(W_{\infty})$ is actually isomorphic to $\mathfrak{sp}_{\infty}$, the Lie algebra of infinite symplectic matrices with finitely non zero entries \cite{Kac94}. This has been discussed in \cite{Z2}, and it plays a key role in proving the simplicity of $V_{\mathcal{J},r}$.
     }
     \par{
        By (\ref{bracket}), we write down the explicit commutation relation for all $a(m)b(n),u(p)v(q)\in S^2(W_{\infty})$, $m,n\in\mathbb{Z}$, which is useful for computations:
        \begin{align}
            &[a(m)b(n),u(p)v(q)]\notag\\
            =&n\delta_{n+p,0}\langle b,u \rangle a(m)v(q)+m\delta_{m+p,0}\langle a,u \rangle b(n)v(q)\notag\\
            +&n\delta_{n+q,0}\langle b,v \rangle u(p)a(m)+m\delta_{m+q,0}\langle a,v \rangle u(p)b(n).\label{bracket3}
        \end{align}
    }
    \par{
        Let $\mathcal{L}(\mathfrak{h})$ denote the following direct sum of Lie algebras
        $$
            \mathcal{L}(\mathfrak{h})\stackrel{def.}{=}S^2(\widetilde{W_{\infty}})\bigoplus \mathbb{C}K,
        $$
        where $K$ is the central element in $\mathcal{L}(\mathfrak{h})$. Define
        \begin{align*}
            &\quad \mathfrak{B}_{+}\stackrel{def.}{=}{\rm span}\{L_{a,b}(m,n)|\,n\geq 0\;\text{or}\; m\geq 0\},\\
            &\mathcal{L}_{-}\stackrel{def.}{=}{\rm span}\{L_{a,b}(m,n)|\,m,n<0\},\;\;\;
            \mathcal{L}_{+}\stackrel{def.}{=}\mathfrak{B}_{+}\bigoplus \mathbb{C}K.
        \end{align*}
        Then we have a decomposition of $\mathcal{L}(\mathfrak{h})$:
        \begin{align*}
        	&\mathcal{L}(\mathfrak{h})=\mathcal{L}_{-}\bigoplus\mathcal{L}_{+}=\mathcal{L}_{-}\bigoplus\mathfrak{B}_{+}\bigoplus \mathbb{C}K.
        \end{align*}
        As mentioned in \cite{Z2}, $\mathcal{L}_{+}$ is certain `parabolic subalgebra' of $\mathcal{L}(\mathfrak{h})$ \cite{KR2}.
    }
    \par{
        Now we consider a map from $\mathcal{L}(\mathfrak{h})$ to $U(\hat{\mathfrak{h}})/(c-1)$. Let $(c-1)$ be the two-sided ideal of $U(\hat{\mathfrak{h}})$ generated by $c-1$. We consider the following map
        $$
            \iota:a(m)b(n)\mapsto \frac{1}{2}(a(m)b(n)+b(n)a(m)),\,S^2(\widetilde{W_{\infty}})\rightarrow U(\hat{\mathfrak{h}}).
        $$
        It is checked that
        $$
            c\iota([x,y])=[\iota(x),\iota(y)]
        $$
        for all $x,y\in S^2(\widetilde{W_{\infty}})$. In particular $\iota$ is a Lie algebra embedding if we take $c=1$, and $\iota$ becomes a map from $S^2(\widetilde{W_{\infty}})$ to $U(\hat{\mathfrak{h}})/(c-1)$. We remark that the `$\frac{1}{2}$' on the left hand side guarantees that $\iota$ is indeed a Lie algebra embedding. Hence the Lie algebra $S^2(\widetilde{W_{\infty}})$ also acts on the Fock space $S(\hat{\mathfrak{h}}_{-})$ by
        $$
            (a(m)b(n))\cdot v\stackrel{def.}{=}\frac{1}{2}(a(m)b(n)+b(n)a(m))v
        $$
        for any $a(m)b(n)\in S^2(\widetilde{W_{\infty}})$ and $v\in S(\hat{\mathfrak{h}}_{-})$. We also note that $\mathcal{L}_{+}\cap \mathfrak{sp}_{\infty}$ is also a parabolic subalgebra of $\mathfrak{sp}_{\infty}$, and
        $\mathbb{C}1$ is a one-dimensional $\mathcal{L}_{+}\cap \mathfrak{sp}_{\infty}$-module. More explicitly we have
        \begin{align}
            (a(m)b(n))\cdot 1=\frac{1}{2}|m|\delta_{m+n,0}(a,b)1\label{spaction}
        \end{align}
        for all $a(m)b(n)\in\mathcal{L}_{+}\cap \mathfrak{sp}_{\infty}$.
    }
    \par{
        Set
        $$
            L_{a,b}(m,n)\stackrel{def.}{=}\frac{1}{2}a(m)b(n)-\frac{1}{4}|m|\delta_{m+n,0}(a,b)K\in \mathcal{L}(\mathfrak{h}).
        $$
        We note that $L_{a,b}(m,n)$ is well defined and
        \begin{align*}
            &L_{a,b}(m,n)=L_{b,a}(n,m),
            &\mathcal{L}(\mathfrak{h})=span\{L_{a,b}(m,n),K\,|\,a,b\in\mathfrak{h},m,n\in\mathbb{Z}\}.
        \end{align*}
        We also note that if $K$ acts as the identity, then the `vacuum vector' in $S(\hat{\mathfrak{h}}_{-})$, then $\mathbb{C}1$ is a one dimensional $\mathcal{L}_{+}$-module such that
        \begin{align}
            L_{a,b}(m,n)\cdot 1=0\text{ for all }m\geq 0\text{ or }n\geq 0. \label{action}
        \end{align}
        The definition of $L_{a,b}(m,n)$ and the Lie algebra $\mathcal{L}=\mathcal{L}(\mathfrak{h})$ are actually the same as the ones given in \cite{Z2} and \cite{AM}, although they looks different.
    }
    \par{
        The key point is that when restricts to $\mathcal{L}_{+}$, we can define a new one dimensional $\mathcal{L}_{+}$-module related to (\ref{action}) such that the central element $K$ acts as $r$ where $r$ is arbitrary complex number. Therefore we are able to define a $\mathcal{L}(\mathfrak{h})$-module $M_r$ which is induced from this one dimensional $\mathcal{L}_{+}$-module. We describe it explicitly as follows. Let $\mathbb{C}v_r$ be the one dimensional $\mathcal{L}_{+}$-module spanned by the element $v_r$ such that
        \begin{align*}
        	&x v_r=0,\text{ for all } x\in\mathfrak{B}_{+},\notag\;\;\;K v_r=rv_r.\notag
        \end{align*}
        It is easy to check that this indeed gives an $\mathcal{L}_{+}$-module. Equivalently, $\mathbb{C}v_r$ is also a one dimensional $\mathfrak{sp}_{\infty}\cap \mathcal{L}_{+}$-module such that
        $$
            (a(m)b(n))\cdot v_r\stackrel{def.}{=}\frac{r}{2}|m|\delta_{m+n,0}(a,b)v_r
        $$
        for all $a(m)b(n)\in\mathfrak{sp}_{\infty}\cap \mathcal{L}_{+}$, which can be compared with (\ref{spaction}). Then we have an $\mathcal{L}(\mathfrak{h})$-module $M_r$ which is induced from this one dimensional $\mathcal{L}_{+}$-module:
        \begin{align}
            M_r\stackrel{def.}{=}&U(\mathcal{L}(\mathfrak{h}))\otimes_{U(\mathcal{L}_{+})}\mathbb{C}v_r
            \cong U(\mathcal{L}_{-})1\notag\\
            =&{\rm span}\{L_{a_1,b_1}(-m_1,-n_1)\cdots L_{a_k,b_k}(-m_k,-n_k)\cdot v_r|\notag\\
            &m_i,n_i\in \mathbb{Z}_{\geq 1},a_i,b_i\in \mathfrak{h}\}.\label{Mr1}
        \end{align}
        We also note that $M_r$ can be viewed as an induced $\mathfrak{sp}_{\infty}$-module
        $$
            M_r=U(\mathfrak{sp}_{\infty})\otimes_{U(\mathfrak{sp}_{\infty})\cap \mathcal{L}_{+}}\mathbb{C}v_r,
        $$
        and this will be used when discussing the simplicity of the VOA.
    }
    \par{
		For $a,b\in\mathfrak{h}$, define the operators $L_{a,b}(l)$ and the fields $L_{a,b}(z)$
		\begin{align}
            L_{a,b}(l)\stackrel{def.}{=}\sum_{k\in\mathbb{Z}}L_{a,b}(-k+l-1,k),\quad
            L_{a,b}(z)\stackrel{def.}{=}\sum_{l\in\mathbb{Z}}L_{a,b}(l)z^{-l-1}.\label{field}
		\end{align}
        It is proved in \cite{AM} that these fields are mutually local. Therefore by reconstruction theorem \cite{Kac}, these mutually local fields generate a vertex algebra:
		$$
			V_{\mathcal{J},r}\stackrel{def.}{=}{\rm span}\{L_{a_1,b_1}(m_1)\cdots L_{a_k,b_k}(m_k)\cdot v_r|\,m_i\in \mathbb{Z},a_i,b_i\in \mathfrak{h}\},
		$$
        and $v_r$ being the `vacuum'. The vertex algebra $V_{\mathcal{J},r}$ is actually a VOA. The Virasoro element $\omega$ is given by:
        $$
            \omega=\sum_k L_{e_k,e_k}(-1,-1)\cdot v_r,
        $$
        and $L(0)=\omega(1)$ gives a gradation on $V_{\mathcal{J},r}$
        $$
            V_{\mathcal{J},r}=\bigoplus_{k\geq 0}(V_{\mathcal{J},r})_{k}.
        $$
        We also check that
        $$
            (V_{\mathcal{J},r})_{0}=\mathbb{C}1, \quad(V_{\mathcal{J},r})_{1}=\{0\},
        $$
        and the Griess algebra $(V_{\mathcal{J},r})_2$ is isomorphic to $\mathcal{J}$:
        $$
            L_{a,b}(-1,-1)\cdot 1\mapsto L_{a,b}\stackrel{def.}{=}a\otimes b+b\otimes a.
        $$
        It was further shown in \cite{NS} that $V_{\mathcal{J},r}=M_r$, but we use $V_{\mathcal{J},r}$ instead of $M_r$ to emphasize the corresponding VOA structure.
    }
    \par{
        We give a set of mutually orthogonal Virasoro elements which corresponds to the Jordan frame. We note that the Greiss algebra $(V_{\mathcal{J},r})_2\simeq S^2(\mathfrak{h})$ is $\frac{d(d+1)}{2}$ dimensional, which is spanned by
        $$
            J(\mathfrak{h})\stackrel{def.}{=}\{L_{e_i,e_i},\,L_{e_i,e_j}|\,1\leq i,j\leq n,\,i\neq j\}.
        $$
        It is easy to check that the following subset
        $$
            F(\mathfrak{h})\stackrel{def.}{=}\{L_{e_i,e_i}|\,1\leq i\leq d\}
        $$
        of $J(\mathfrak{h})$ consists of mutually orthogonal Virasoro elements by a direct computation. Each element in $F(\mathfrak{h})$ is a Virasoro element of central charge $1$, and the half of these elements are bijective to the the elements in the Jordan frame through the isomorphism between $\mathcal{J}(W)$ and the Greiss algebra $(V_{\mathcal{J},r})_2$:
        $$
            L_{e_i,e_i}\mapsto \frac{1}{2}L_{e_i,e_i},\,\mathcal{J}_B(\mathfrak{h}) \rightarrow (V_{\mathcal{J}},r)_2.
        $$
        We remark that the normalized factor `$\frac{1}{2}$' is used to guarantee that this is an algebra homomorphism.
    }
    \par{
        Let $\bar{V}_{\mathcal{J},1}$ be the VOA constructed in \cite{Lam1} which is obtained by taking the $-1$ fixpoint of $S(\hat{\mathfrak{h}}_{-})$, where `$-1$' denotes the action induced by multiplying $-1$ on $\mathfrak{h}$. In \cite{Z2} we showed that $\bar{V}_{\mathcal{J},1}$ is the simple quotient of $V_{\mathcal{J},1}$. We note that the properties of $V_{\mathcal{J},r}$ are quite similar to $\bar{V}_{\mathcal{J},1}$, therefore we can study $V_{\mathcal{J},r}$ by comparing it with $\bar{V}_{\mathcal{J},1}$. For example, for $a(-1)b\in \bar{V}_{\mathcal{J},1}\subseteq S(\hat{\mathfrak{h}}_{-})$, the corresponding field is given by
        $$
            Y(a(-1)b,z)=\sum_{l\in \mathbb{Z}} (a(-1)b)(l)z^{-l-1}=\sum_{k,l\in \mathbb{Z}} :a(-k+l-1)b(k):z^{-l-1}.
        $$
        If we replace $a(-1)b$ with $2L_{a,b}$ and $:a(m)b(n):$ with $2L_{a,b}(m,n)$ in the above, it is exactly the same as (\ref{field}). We can also get some other identities starting from the identities in $V_{\mathcal{J},1}$ (and therefore in $S(\hat{\mathfrak{h}}_{-})$). The difference appears when the central element $c$ arises in the computation. For example, for $a(-1)b,u(-1)v\in \bar{V}_{\mathcal{J},1}$, we have
        \begin{align*}
            (a(-1)b,u(-1)v)=(a(-1)b)(1)u(-1)v=(a,u)(b,v)+(a,v)(b,u),
        \end{align*}
        while for $2L_{a,b},2L_{u,v}\in V_{\mathcal{J},r}$, we have
        \begin{align*}
            (2L_{a,b},2L_{u,v})=(2L_{a,b})(1)2L_{u,v}=r(a,u)(b,v)+r(a,v)(b,u)=r(a(-1)b,u(-1)v).
        \end{align*}
        In \cite{Z1} we computed the correlation functions of generating fields in $V_{\mathcal{J},r}$. If we view the correlation function as a function of $r$, the result is a polynomial function of $r$.
    }
    \section{Construction of the VOA $V_{\mathcal{J},r}$, Where $\mathcal{J}$ is a Type $C$ Jordan Algebra}
    \par{
        Let $r$ be arbitrary complex number. In this section, we construct the VOA $V_{\mathcal{J},r}$ satisfying the property that $(V_{\mathcal{J},r})_0=\mathbb{C}1$, $(V_{\mathcal{J},r})_1=\{0\}$, and $(V_{\mathcal{J},r})_2$ is isomorphic to the type $C$ Jordan algebra $\mathcal{J}_C(W)$.
    }
    \par{
        Recall that in Section 3, we realize the type $B$ simple Jordan algebra $\mathcal{J}_B(\mathfrak{h})$ using the symmetric square space $S^2(\mathfrak{h})$, and the type $C$ simple Jordan algebra $\mathcal{J}_C(W)$ is realized in a similar way using the symmetric wedge $\bigwedge^2(W)$. The analogy between these two types of Jordan algebras indicates the construction of $V_{\mathcal{J},r}$ for $\mathcal{J}\simeq \mathcal{J}_C(W)$, which is almost parallel to Section 3.
    }
    \par{
        In this section, we fix the notation that $W$ denotes a $2n$-dimensional symplectic space with the symplectic form $\langle\cdot,\cdot\rangle$, and $\mathcal{J}\simeq \mathcal{J}_C(W)$. We cite the following formulas for a given VOA $V$ without proof:
        \begin{align}
            &[x(m),y(n)]=\sum_{j\geq 0}{m\choose j}(x(j)y)(m+n-j), \label{form1}\\
            &(x(n)y)(l)\notag\\
            =&\sum_{i\geq 0}{n \choose i}(-1)^i(x(-i+n)y(i+l)-(-1)^{n}y(-i+n+l)x(i)),\label{form2}\\
            &\text{ for all }x,y\in V \text{ and } m,n,l\in\mathbb{Z}.
        \end{align}
        The details can be found, for example, in \cite{Kac}. These formulas will be frequently used in this section. Here we note that the combinatorial coefficients means
        $$
            {m\choose i}\stackrel{def.}{=}\frac{m\times\cdots (m-i+1)}{i\times\cdots 1}
        $$
        for all $m\in \mathbb{Z},i\in \mathbb{Z}_{\geq 1}$, and
        $$
            {m\choose 0}=1
        $$
        for all $m\in \mathbb{Z}$ by convention. An elementary calculation shows that the recursion relation
        \begin{align}
            {m\choose i}={m-1\choose i-1}+{m-1\choose i} \label{comb}
        \end{align}
        holds for all $m\in \mathbb{Z},i\in \mathbb{Z}_{\geq 1}$.
    }
    \par{
        Given a vector space $V$ with a symmetric bilinear form $(\cdot,\cdot)$, the space $\bigwedge^2(V)$ is a Lie algebra such that for $ab,uv\in \bigwedge^2(W)$, the Lie bracket is given by
        \begin{align}
            [ab,uv]=(b,u) av-(a,u) bv+(b,v) ua-(a,v)ub. \label{bracket2}
        \end{align}
        In particular when $V$ is finite dimensional and $(\cdot,\cdot)$ is non-degenerate, $\bigwedge^2(V)$ is isomorphic to the finite dimensional simple Lie algebra $\mathfrak{so}(V)$. It is well known that this is a special case of the spinor construction of irreducible $\mathfrak{so}(V)$-modules.
    }
    \par{
        We apply this to
        $$
            \widetilde{V_{\infty}}{=}W\otimes\mathbb{C}[t,t^{-1}].
        $$
        Recall the anti-symmetric bilinear form $\langle\cdot,\cdot\rangle'$ over the space of Laurant polynomials $\mathbb{C}[t,t^{-1}]$ given in Section 3. It is obvious that $\widetilde{V_{\infty}}$ has a symmetric bilinear form $(\cdot,\cdot)=\langle\cdot,\cdot\rangle\otimes \langle\cdot,\cdot\rangle'$:
        $$
            ( a(m),b(n))=m\langle a,b\rangle \delta_{m+n,0},
        $$
        hence $\bigwedge^2(\widetilde{V_{\infty}})$ is a Lie algebra with the bracket given by (\ref{bracket2}). The following explicit formula for the commutation relation, which can be easily derived from (\ref{bracket2}), is useful for later computations:
        \begin{align}
            &[a(m)b(n),u(p)v(q)]\notag\\
            =&n\delta_{n+p,0}\langle b,u \rangle a(m)v(q)-m\delta_{m+p,0}\langle a,u \rangle b(n)v(q)\notag\\
            -&n\delta_{n+q,0}\langle b,v \rangle u(p)a(m)+m\delta_{m+q,0}\langle a,v \rangle u(p)b(n),\label{bracket4}\\
            &\text{ for all } a(m)b(n)\in \bigwedge^2(\widetilde{V_{\infty}}),\,m,n\in\mathbb{Z}.\notag
        \end{align}
        We also note that the symmetric bilinear form $(\cdot,\cdot)$ restricted to the subspace
        $$
            V_{\infty}\stackrel{def.}{=}W\otimes\mathbb{C}t[t]\bigoplus W\otimes\mathbb{C}t^{-1}[t^{-1}]
        $$
        is non-degenerate, and the Lie subalgebra $\bigwedge^2(V_{\infty})$ is actually isomorphic to $\mathfrak{so}_{\infty}$, the Lie algebra of infinite orthogonal matrices \cite{Kac94}.
    }
    \par{
        We also recall some basic facts about an infinite dimensional Lie superalgebra associated to $W$, and the corresponding SVOA called the symplectic Fermion SVOA. Note that there is a Lie superalgebra $\widehat{W}$ associated to $W$:
        \begin{align*}
            \widehat{W}=\widetilde{V_{\infty}}\oplus \mathbb{C}c.
        \end{align*}
        The even part $\widehat{W}_{\bar{0}}$ and the odd part $\widehat{W}_{\bar{1}}$ are
        $$
            \widehat{W}_{\bar{0}}=\mathbb{C}c,\,\widehat{W}_{\bar{1}}=\widetilde{V_{\infty}},
        $$
        and the Lie superbracket is given by
        \begin{align*}
            [a(m),b(n)]=(a(m),b(n))c=m\langle a,b\rangle \delta_{m+n,0}c,\;\;
            [x,c]=0,\text{ for all } x\in \widehat{W},\notag
        \end{align*}
        where $a(m)\stackrel{def.}{=}a\otimes t^m$, $a\in W$. It follows that
        $$
            \widehat{W}_{-}\stackrel{def.}{=}W \otimes \mathbb{C}t^{-1}[t^{-1}]
        $$
        is a super-commutative Lie sub-superalgebra of $\widehat{W}$. The enveloping algebra $U(\widehat{W})$ is essentially certain kind of infinite dimensional Clifford algebra, and the Fock space $\bigwedge(\widehat{W}_{-})\cdot 1\simeq \bigwedge(\widehat{W}_{-})$ is a left $U(\widehat{W})$-module. The space $\bigwedge(\widehat{W}_{-})$ has a SVOA structure, and the field associated to $x=x(-1)\cdot 1$ is:
        $$
            Y(x,z)=\sum_{k}x(k)z^{-k-1},\quad x\in W.
        $$
        This is called the `symplectic Fermion' SVOA in literatures \cite{Ab07}, denoted by
        $\mathcal{A}(W)$.
    }
    \par{
        Let
        $$
            \mathcal{L}(W)\stackrel{def.}{=}\bigwedge^2(\widetilde{V_{\infty}})\bigoplus \mathbb{C}K,
        $$
        where $K$ is the central element:
        $$
            [x,K]=0,\text{ for all }x\in \bigwedge^2(\widetilde{V_{\infty}}).
        $$
        We note that there is a map
        $$
            \iota:a(m)b(n)\mapsto \frac{1}{2}(a(m)b(n)-b(n)a(m))
        $$
        from $\bigwedge^2(\widetilde{V_{\infty}})$ to $U(\widehat{W})$. It is easy to check that $\iota$ satisfies
        \begin{align}
            c\iota([x,y])=[\iota(x),\iota(y)].\label{deform2}
        \end{align}
        In particular, if we view $\iota$ as a map from $\bigwedge^2(\widetilde{V_{\infty}})$ to $U(\widehat{W})/(c-1)$, then $\iota$ is a Lie algebra embedding. Therefore $\bigwedge^2(\widetilde{V_{\infty}})$ acts on the Fermionic Fock space $\bigwedge(\widehat{W}_{-})\cdot 1$.
    }
    \par{
        Set
        \begin{align}
            L_{a,b}(m,n)\stackrel{def.}{=}a(m)b(n)-\frac{1}{2}|m|\delta_{m+n,0}\langle a,b\rangle K\in \mathcal{L}(W).\label{normal}
        \end{align}
        The element $L_{a,b}(m,n)$ is also well defined and
        \begin{align*}
            &L_{a,b}(m,n)=-L_{b,a}(n,m),&\mathcal{L}(W)=span\{L_{a,b}(m,n),K\,|\,a,b\in W,m,n\in\mathbb{Z}\}.
        \end{align*}
        Note that there is a sign change for $L_{a,b}(m,n)$ when we exchange the pairs $(a,m),(b,n)\in W\times \mathbb{Z}$. Define
        \begin{align*}
            &\quad \mathfrak{B}_{+}\stackrel{def.}{=}{\rm span}\{L_{a,b}(m,n)|\,n\geq 0\;\text{or}\; m\geq 0\},\\
            &\mathcal{L}_{-}\stackrel{def.}{=}{\rm span}\{L_{a,b}(m,n)|\,m,n<0\},\;\;\;
            \mathcal{L}_{+}\stackrel{def.}{=}\mathfrak{B}_{+}\bigoplus \mathbb{C}K.
        \end{align*}
        Then we have a decomposition of $\mathcal{L}(W)$:
        \begin{align*}
        	&\mathcal{L}(W)=\mathcal{L}_{-}\bigoplus\mathcal{L}_{+}=\mathcal{L}_{-}\bigoplus\mathfrak{B}_{+}\bigoplus \mathbb{C}K.
        \end{align*}
    }
    \par{
        Similar to what we have done in the type $B$ case, we can also construct a one dimensional $\mathcal{L}_{+}$-module such that $K$ acts as arbitrary complex number $r$. Let $\mathbb{C}v_r$ be the one dimensional $\mathcal{L}_{+}$-module such that
        \begin{align}
        	x v_r=0,\text{ for all } x\in\mathfrak{B}_{+},Kv_r=rv_r.\label{action1}
        \end{align}
        Then we construct $M_r$ as an induced $\mathcal{L}(W)$-module:
        \begin{align}
            M_r\stackrel{def.}{=}&U(\mathcal{L}(W))\otimes_{U(\mathcal{L}_{+})}\mathbb{C}v_r
            \cong U(\mathcal{L}_{-})1\notag\\
            =&{\rm span}\{L_{a_1,b_1}(-m_1,-n_1)\cdots L_{a_k,b_k}(-m_k,-n_k)\cdot v_r|\notag\\
            &m_i,n_i\in \mathbb{Z}_{\geq 1},a_i,b_i\in W\}.\label{Mr1}
        \end{align}
        We note that $\mathbb{C}v_r$ is also a one dimensional $\mathfrak{so}_{\infty}\cap\mathcal{L}_{+}$-module. The action is given by
        \begin{align}
            (a(m)b(n))\cdot v_r=\frac{r}{2}|m|\delta_{m+n,0}\langle a,b\rangle v_r \label{action2}
        \end{align}
        for all $a(m)b(n)\in \mathfrak{so}_{\infty}\cap \mathcal{L}_{+}$ by (\ref{normal}), and $M_r$ is also an induced $\mathfrak{so}_{\infty}$-module:
        $$
            M_r=U(\mathfrak{so}_{\infty})\otimes_{U(\mathfrak{so}_{\infty})\cap \mathcal{L}_{+}}\mathbb{C}v_r.
        $$
        This will be used later in Section 6.
    }
    \par{
		For $a,b\in W$, define the operators $L_{a,b}(l)$ and the fields $L_{a,b}(z)$
		\begin{align}
            L_{a,b}(l)\stackrel{def.}{=}\sum_{k\in\mathbb{Z}}L_{a,b}(-k+l-1,k),\quad
            L_{a,b}(z)\stackrel{def.}{=}\sum_{l\in\mathbb{Z}}L_{a,b}(l)z^{-l-1}.\label{mode}
		\end{align}
        Then we have
        \begin{proposition}
            The formal power series $L_{a,b}(z),\,a,b\in W$ are mutually local fields over $M_r$.
        \end{proposition}
    }
    \par{
        \textbf{Proof of Proposition 1.} The proof is essentially the same as the proof given in \cite{AM} for type $B$ case. First it easy to show that for a given $v\in M_r$, there exists an integer $N$ (depending on $a$, $b$ and $v$) such that
        $$
            L_{a,b}(n)v=0
        $$
        for all $n\geq N$, hence $L_{a,b}(z)$ are all fields. We also check that for all $n\geq 0$,
        $$
            L_{a,b}(n)\cdot 1=0.
        $$
    }
    \par{
        Next we show the locality by comparing the fields $L_{a,b}(z)$ with the fields $Y(a(-1)b,z)$ in the symplectic Fermion SVOA $\mathcal{A}(W)$. It is easy to check that for all $a,b\in W$ the $U(\widehat{W})$-valued formal power series are mutually local:
        $$
            (z-w)^2[a(z),b(w)]=0.
        $$
        By a variant of Dong's lemma for Lie superalgebra valued formal power series \cite{Kac}, for any $a,b,u,v\in W$, the $U(\widehat{W})$-valued formal power series $Y(a(-1)b,z)$ and $Y(u(-1)v,z)$ are also mutually local:
        $$
            (z-w)^4[Y(a(-1)b,z),Y(u(-1)v,z)]=0.
        $$
    }
    \par{
        By (\ref{deform2}), (\ref{normal}), and (\ref{mode}) we have
        $$
            c(z-w)^4\iota([L_{a,b}(z),L_{u,v}(w)])=[\iota(Y(a(-1)b,z)),\iota(Y(u(-1)v,w))]=0.
        $$
        as $U(\widehat{W})$-valued formal power series. Because $\iota$ is an embedding and $cx=0$ implies $x=0$ in $U(\widehat{W})$, therefore we have
        $$
            (z-w)^4[L_{a,b}(z),L_{u,v}(w)]=0,
        $$
        and we conclude the proof of the locality of the fields. $\qed$
    }
    \par{
        By Proposition 1 and the reconstruction theorem \cite{Kac}, these mutually local fields generate a vertex algebra:
		$$
			V_{\mathcal{J},r}\stackrel{def.}{=}{\rm span}\{L_{a_1,b_1}(m_1)\cdots L_{a_k,b_k}(m_k)\cdot 1|\,m_i\in \mathbb{Z},a_i,b_i\in \mathfrak{h}\}.
		$$
    }
    \par{
        Now we check the remaining part of (1) in Theorem 1 about the VOA $V_{\mathcal{J},r}$. We note that to do the computations, we only need (\ref{bracket4}) and (\ref{normal}). For any $a,b,u,v\in W$, a direct computation shows that
        $$
            (L_{a,b})(1)(L_{u,v})=\langle b,u\rangle L_{a,v}-\langle b,v\rangle L_{a,u}-\langle a,u\rangle L_{b,v}+\langle a,v\rangle L_{b,u}.
        $$
        Compared with (\ref{JordanCmul}) it is clear that
        $$
            L_{a,b}\mapsto L_{a,b},\,(V_{\mathcal{J},r})_2\mapsto \mathcal{J}(W)
        $$
        is an algebra isomorphism.
    }
    \par{
        Now we give a set of mutually orthogonal Virasoro elements in $(V_{\mathcal{J},r})_2$, which corresponds to the elements in the Jordan frame of $\mathcal{J}(W)$ mentioned in Section 2. We note that the Greiss algebra $(V_{\mathcal{J},r})_2\simeq \bigwedge^2(W)$ is $n(2n-1)$ dimensional, which is spanned by
        $$
            J(W)\stackrel{def.}{=}\{L_{\psi_i,\psi^{*}_i},\,L_{\psi_i,\psi_j},\,L_{\psi^{*}_i,\psi^{*}_j},\,L_{\psi_j,\psi^{*}_k}|\,1\leq i,j,k\leq n,\,j\neq k\}.
        $$
        It is easy to check that the following subset
        $$
            F(W)\stackrel{def.}{=}\{L_{\psi_i,\psi^{*}_i}|\,1\leq i\leq n\}
        $$
        of $J(W)$ consists of mutually orthogonal Virasoro elements by a direct computation. Each element in $F(W)$ is a Virasoro element of central charge $-2$, and the half of these elements are bijective the the elements in the Jordan frame through the isomorphism between $\mathcal{J}(W)$ and the Greiss algebra $V_{\mathcal{J},r}$:
        $$
            \frac{1}{2}L_{\psi_i,\psi^{*}_i}\mapsto L_{\psi_i,\psi^{*}_i},\,(V_{\mathcal{J},r})_2\rightarrow \mathcal{J}_C(W).
        $$
        Moreover, the sum of these elements
        $$
            \omega\stackrel{def.}{=}\sum_i L_{\psi_i,\psi^{*}_i}
        $$
        is also a Virasoro element, with central charge equals $-2n$. The element $\omega(1)$ gives the gradation of $V_{\mathcal{J},r}$:
        $$
            V_{\mathcal{J},r}=\sum_{k\geq 0}(V_{\mathcal{J},r})_k,
        $$
        where $k$ is the eigenvalue of $\omega(1)$. It follows easily from the construction of $V_{\mathcal{J},r}$ that $(V_{\mathcal{J},r})_0=\mathbb{C}1$ and $(V_{\mathcal{J},r})_1=\{0\}$. hence the Greiss algebra $(V_{\mathcal{J},r})_2$ is indeed isomorphic to the type $C$ Jordan algebra $\mathcal{J}(W)$, and we conclude the proof of (1) in Theorem 1.
    }
    \par{
        Let $\bar{V}_{\mathcal{J},1}$ be the sub VOA of $\bigwedge(\hat{W})_{-}$ by taking the $-1$ fixpoint of $\bigwedge(\hat{W})_{-}$, where `$-1$' denotes the action induced by multiplying $-1$ on $W$. It follows that $V_{\mathcal{J},1}$ acts on $\bar{V}_{\mathcal{J},1}$ and $\bigwedge(\hat{W})_{-}$ such that
        $$
            L_{a,b}(m,n)=:a(m)b(n):,\,K=Id
        $$
        for all $a,b\in W,m,n\in\mathbb{Z}$, and $\bar{V}_{\mathcal{J},1}$ is essentially a quotient of $\bar{V}_{\mathcal{J},1}$. The VOA $\bar{V}_{\mathcal{J},1}$ is studied by Abe in \cite{Ab07}, and it will also plays a role when discussing the simplicity of $V_{\mathcal{J},r}$.
    }
    \par{
        We give some formulas to conclude this section, which are needed later. We note that for $i\neq k$, $L_{\psi_k,\psi^{*}_k}$ and $L_{\psi_l,\psi^{*}_l}$ are two mutually orthogonal Virasoro elements in $F(W)$, $L_{\psi_k,\psi^{*}_l}\in J(W)$ is not an element in $F(W)$. By (\ref{form1}) we have
        $$
            (x(0)y)(l)=[x(0),y(l)]
        $$
        for any $x,y$ in a VOA $V$ and $l\in\mathbb{Z}$, hence it is easy to show that
        \begin{align}
            &(L_{\psi_k,\psi^{*}_l}(-i-1,-j)1)(l)=\frac{1}{i}[L_{\psi_k,\psi^{*}_l}(0),(L_{\psi_k,\psi^{*}_l}(-i,-j)1)(l)],\label{vircom1}\\
            &(L_{\psi_k,\psi^{*}_l}(-i,-j-1)1)(l)=\frac{1}{j}[L_{\psi_l,\psi^{*}_l}(0),(L_{\psi_k,\psi^{*}_l}(-i,-j)1)(l)]\label{vircom2}
        \end{align}
        for all $i,j\geq 1$. Apply this repeatedly together with (\ref{form1}), (\ref{form2}) and (\ref{mode}),
        we have the following identities:
        \begin{lemma}
            If $L_{a,b}\in J(W)$ is not an element in $F(W)$, then for all $s\in\mathbb{Z}$ and $i,j\geq 1$ the following formula holds
            \begin{align}
                &(L_{a,b}(-i,-j)1)(s+i+j-1)\notag\\
                =&\sum_{k\in \mathbb{Z}}(-1)^{i+j}{i+k-1\choose i-1}{j+s-k-1\choose j-1}L_{a,b}(k,s-k).\label{nonvir0}
            \end{align}
            In particular when $j=1$,
            \begin{align}
                (L_{a,b}(-i,-1)1)(s+i)=\sum_{k\in \mathbb{Z}}(-1)^{i-1}{i+k-1\choose i-1}L_{a,b}(k,s-k).\label{nonvir}
            \end{align}
        \end{lemma}
        We remark that (\ref{nonvir0}) and (\ref{nonvir}) actually holds for any $L_{a,b}\in J(W)$, but we do not need this fact later.
    }
    \section{Properties of the VOA $V_{\mathcal{J},r}$, Where $\mathcal{J}$ is a Type $C$ Jordan Algebra}
    \par{
        In this section, we study the property of the VOA $V_{\mathcal{J},r}$, where $\mathcal{J}=\mathcal{J}_C(W)$. We will prove that $V_{\mathcal{J},r}$-submodules of $M_r$ are essentially the same as $\mathcal{L}(W)$-submodules of $M_r$. This result is analogous to Proposition 3.4 in \cite{NS}. Therefore, the study of the $V_{\mathcal{J},r}$ can be reduced to the study of $M_r=V_{\mathcal{J},r}$ as an $\mathcal{L}(W)$-module. We will also show that (2) of Theorem 1 will be a corollary of Proposition 1 at the end of this section.
    }
    \par{
        The proofs in this section are similar to the one in \cite{NS}, and we just make some modifications. We first need the following lemma:
        \begin{lemma}
            For any $m,n\in \mathbb{Z}$ and $L_{a,b}\in J(W)$, but $L_{a,b}\notin F(W)$, we have
            $$
                L_{a,b}(m,n)u\in span\{L_{a,b}(-p,-1)(l)u|\,a,b\in W,\,p\in\mathbb{Z}_{\geq 1},l\in\mathbb{Z}\}.
            $$
        \end{lemma}
    }
    \par{
        \textbf{Proof of Lemma 2.} Because $L_{a,b}(m,n)$ acts as zero on $M_r$ if $m=0$ or $n=0$, so we only need to consider the case $m,n\neq 0$. Because for a given element $u$, there exists an integer $t$ and a positive integer $N$ such that
        $$
            L_{a,b}(j,s-j)u=0
        $$
        for all $s<t$ and $s>t+N$, therefore by (\ref{nonvir})
        $$
            (L_{a,b}(-i,-1)1)(s+i)u=\sum_{t\leq j\leq t+N}(-1)^{i-1}{i+j-1\choose i-1}L_{a,b}(j,s-j)u.
        $$
        That is, $(L_{a,b}(-i,-1)1)(s+i)u$ is a linear combination of $L_{a,b}(j,s-j)u$.
    }
    \par{
        Now we show that we can express $L_{a,b}(j,s-j)u,t\leq j\leq t+N$ in terms of $(L_{a,b}(-i,-1)1)(s+i)u,i=1,\cdots N+1$. It is enough to show that the $(N+1)\times (N+1)$ matrix $(a^{T,N}_{ij}),1\leq i,j\leq N+1$ where
        $$
            a^{t,N}_{i,j}=(-1)^{i-1}{t+i+j-2\choose i-1}
        $$
        has non-zero determinant. This can be shown by elementary linear algebra as follows. Let $(a'_{i,j}),1\leq i,j\leq N+1$ be the new matrix such that
        $$
            a'_{i,j}=a^{t,N}_{i,j}-a^{t,N}_{i,j-1},a'_{i,1}=a_{i,1}
        $$
        for all $2\leq j\leq N+1$, that is, the $j$-th column of $(a'_{i,j})$ equals the $j$-th column of $(a^{t,N}_{i,j})$ subtracting the $j-1$-th column of $(a^{t,N}_{i,j})$, $2\leq j\leq N+1$, and the first column of $(a'_{i,j})$ equals the first column of $(a^{t,N}_{i,j})$. It follows from (\ref{comb}) that
        $$
             a'_{i,j}=-a^{t,N}_{i-1,j}
        $$
        for all $2\leq i,j\leq N+1$, and
        $$
            a'_{1,1}=1,\,a'_{1,j}=0
        $$
        for all $2\leq j \leq N+1$. By the property that a determinant of a given matrix is invariant under elementary operations, we see that
        $$
            det((a^{t,N}_{i,j}))=det(a'_{i,j})=(-1)^Ndet(a^{t+1,N-1}_{i,j}).
        $$
        It is also clear that when $N=0$,
        $$
             det((a^{t,0}_{i,j}))=1
        $$
        for all $t\in\mathbb{Z}$, hence
        $$
            det((a^{t,N}_{i,j}))=(-1)^{\frac{N(N+1)}{2}}\neq 0,
        $$
        and we conclude the proof. $\qed$.
        }
        \par{
        We also need the following lemma which can be derived from the previous two lemmas:
        \begin{lemma}
            For any $a,b\in W$, $m,n\in \mathbb{Z}$ and $u\in M_r$, we have
            \begin{align*}
                &span\{L_{a,b}(m,n)u|\,a,b\in W,m,n\in\mathbb{Z}_{\geq 1}\}\\
                \subseteq &span\{L_{a_1,b_1}(l_1)\cdots L_{a_k,b_k}(l_k)u|\,a_i,b_i\in W,\,l_i\in\mathbb{Z}\}.
            \end{align*}
        \end{lemma}
        \textbf{Proof of Lemma 3.} We first prove the case for $L_{a,b}\in J(W)$ but $L_{a,b}\notin F(W)$.
        It follows from (\ref{vircom1}) and (\ref{vircom2}) that for all $m,n\geq 1$ and $l\in\mathbb{Z}$ we have
        \begin{align}
            (L_{a,b}(-m,-n)1)(l)
            \in span\{L_{a_1,b_1}(l_1)\cdots L_{a_k,b_k}(l_k)|\,a_i,b_i\in W,\,l_i\in\mathbb{Z}\}
        \end{align}
        using induction on $m$ and $n$. Then by Lemma 2, it is clear that
        \begin{align}
            L_{a,b}(m,n)u \in &span\{L_{a,b}(-p,-1)(l)u|\,a,b\in W,\,p\in\mathbb{Z}_{\geq 1},l\in\mathbb{Z}\}\notag\\
            \subseteq&span\{L_{a_1,b_1}(l_1)\cdots L_{a_k,b_k}(l_k)u|\,a_i,b_i\in W,\,l_i\in\mathbb{Z}\}, \label{belong}
        \end{align}
        for all $m,n\in\mathbb{Z}$, so we finish the proof for this case.}
        \par{
             Now we consider the case $L_{a,b}\in F(W)$. It is enough to show that
         $$
            L_{\psi_{1},\psi^{*}_1}(m,n)u
            \in span\{L_{a_1,b_1}(l_1)\cdots L_{a_k,b_k}(l_k)u|\,a_i,b_i\in W,\,l_i\in\mathbb{Z}\}.
            $$
        for all $m,n\in\mathbb{Z}$. This follows easily from the fact that
        $$
            L_{\psi_{1},\psi^{*}_1}(m,n)=[L_{\psi_{1},\psi^{*}_2}(m,1),L_{\psi_{2},\psi^{*}_1}(-1,n)]
        $$
        holds for all $m,n\in\mathbb{Z}$, together with (\ref{belong}) shown above. Hence we conclude the proof of Lemma 3. $\qed$
    }
    \par{
        It can be seen from the proof that the condition $dim(W)\geq 4$ is necessary, because when $dim(W)=2$, all elements in $J(W)$ are Virasoro elements, hence we cannot find non Virasoro elements in $J(W)$ to guarantee the first half in the proof of Lemma 3. Some discussion about the case $dim(W)=2$ are given, for example, in \cite{Ab07}.
    }
    \par{
        The next proposition is an analogue of the Proposition 3.4 in \cite{NS}, which will be used in next section.
        \begin{proposition}
            Suppose $dim(W)\geq 4$, then any $\mathcal{L}(W)$-submodule of $M_r$ is also a $V_{\mathcal{J},r}$-module. Conversely, any $V_{\mathcal{J},r}$-submodule of $M_r$ is a $\mathcal{L}(W)$-module.
        \end{proposition}
    }
    \par{
        \textbf{Proof of Proposition 1.} First, we note that by (\ref{mode}), the `mode' $L_{a,b}(l)$ of the field $Y(L_{a,b},z)$ is an (infinite) sum of elements in $\mathcal{L}(W)$. Suppose $M$ is a $\mathcal{L}(W)$-submodule of $M_r$, we have shown that for any $l\in\mathbb{Z}$ and $a,b\in W$, $L_{a,b}(l)M\subseteq M$ holds, because
        $$
            L_{a,b}(l)u=\sum_{|k|\leq N}L_{a,b}(-k+l-1,k)u\subseteq M
        $$
        for some integer $N$ depending on $L_{a,b}(l)$ and $u$. Therefore any $\mathcal{L}(W)$-submodule of $M_r$ is also a $V_{\mathcal{J},r}$-module. The converse follows from Lemma 3. $\qed$.
    }
    \par{
        \textbf{Proof of (2) in Theorem 1:} Using induction on $k$, it follows easily from Lemma 3 that for any $k\geq 1$,
        \begin{align}
                &span\{L_{a_1,b_1}(-m_1,-n_1)\cdots L_{a_k,b_k}(-m_k,-n_k)1|\,a_i,b_i\in W,m_i,n_i\in\mathbb{Z}_{\geq 1}\}\notag\\
                \subseteq &span\{L_{a_1,b_1}(l_1)\cdots L_{a_k,b_k}(l_k)1|\,a_i,b_i\in W,\,l_i\in\mathbb{Z}\}, \label{span2}
        \end{align}
        which is exactly
        $$
            M_r=V_{\mathcal{J},r}.
        $$
        Therefore we conclude the proof of (2) in Theorem 1. $\qed$
    }
    \par{
        We remark that the analogue of Lemma 3 proved in this section can also be applied to simplify the proofs in \cite{NS} for the case $\mathcal{J}=\mathcal{J}_B(\mathfrak{h})$.
    }
    \section{Simplicity of the VOA $V_{\mathcal{J},r}$, Where $\mathcal{J}$ is a Type $C$ Jordan Algebra}
    \par{
        In this section we prove (3) of Theorem 1 about the simplicity of the VOA $V_{\mathcal{J},r}$, where $\mathcal{J}=\mathcal{J}(W)$. The method we use here is essentially the same as the approach given in \cite{Z2}.
    }
    \par{
        Recall that the infinite dimensional Lie algebra $\bigwedge^2(\widetilde{V_{\infty}})$ (and therefore, $\mathcal{L}(W)$) contains the Lie subalgebra $\bigwedge^2(V_{\infty})$, which is isomorphic to $\mathfrak{so}_{\infty}$. We note that
        $$
            \mathcal{I}\simeq span\{L_{a,b}(0,m)|\,a,b\in W,m\in\mathbb{Z}\}
        $$
        acts as zero on $M_r$, and
        $$
            \mathcal{L}(W)/\mathcal{I}\simeq \mathfrak{so}_{\infty}\bigoplus \mathbb{C}c.
        $$
        Hence, the study of $M_r$ as an $\mathcal{L}(W)$-module is the same as the study of $M_r$ as a $\mathfrak{so}_{\infty}$-module. Combine this with Proposition 1, we have
        \begin{proposition}
            $V_{\mathcal{J},r}$ is a simple VOA if and only if $M_r=V_{\mathcal{J},r}$ is simple as a $\mathfrak{so}_{\infty}$- module.
        \end{proposition}
    }
    \par{
        We first analyze the Lie algebra structure of $\mathfrak{so}_{\infty}$. It is easy to choose countable number of vectors $\{e_i|\,i\in\mathbb{Z}_{\neq 0}\}$ with $(e_i,e_j)=\delta_{i+j,0}$ such that they span $V_{\infty}$, and for all $i\neq 0$, there are some $a,b\in W$ and $j\in\mathbb{Z}_{\geq 1}$ such that
        $$
            e_i=a(j),e_{-i}=b(-j),\,\langle a,b\rangle=1.
        $$

        Then we see that
        $V_{\infty}$ has an increasing chain of subspaces $V_N$:
        $$
            \{0\}=V_0\subseteq V_1\subseteq \cdots V_N\subseteq V_{\infty},\,V_{N}=span\{e_i|\,|i\leq N|\}.
        $$
        It is clear that the symmetric bilinear form restricted to $V_N$ is non-degenerate, and $dim(V_N)=2N$. Because $\bigwedge^2(V_N)\simeq \mathfrak{so}_{2N}$, there is an increasing chain of Lie algebras:
        $$
            \mathfrak{so}_{2}\subseteq \cdots \mathfrak{so}_{2N}\subseteq \mathfrak{so}_{\infty},
        $$
    }
    \par{
        Set
        $$
        \mathfrak{g}^{(N)}\stackrel{def.}{=}\mathfrak{so}(2N)\simeq {\rm span}\{e_k e_l|\,1\leq |k|,|l|\leq N,k\neq l\},
        $$
        where we write
        $$
            e_ke_l\stackrel{def.}{=}e_k\wedge e_l\in \mathfrak{so}_{2N}
        $$
        for convenience. We note that
        $$
        \mathfrak{g}^{(N)}=\mathfrak{g}^{(N)}_{+}\bigoplus \mathfrak{h}^{(N)}\bigoplus \mathfrak{g}^{(N)}_{-},
        $$
        where
        \begin{align*}
        &\mathfrak{g}^{(N)}_{+}={\rm span}\{e_k e_l|\,k+l>0\},\mathfrak{g}^{(N)}_{-}={\rm span}\{e_k e_l|\,k+l<0\},\\
        &\mathfrak{h}^{(N)}= {\rm span}\{e_k e_l|\,k+l=0\}.
        \end{align*}
        Introduce elements $\epsilon_k\in(\mathfrak{h}^{(N)})^{*},k=1,\cdots, N$ such that:
        $$
            \epsilon_l(e_{-k}e_{k})=-\delta_{k,l}.
        $$
        The positive and negative roots with respect to the triangular decomposition are:
        \begin{align*}
        &\Phi^{(N)}_{+}=\{+\epsilon_i+\epsilon_j|\,i<j\}\cup\{-\epsilon_i+\epsilon_j|\,i<j\},\\
        &\Phi^{(N)}_{-}=\{-\epsilon_i-\epsilon_j|\,i<j\}\cup\{+\epsilon_i-\epsilon_j|\,i<j\}.
        \end{align*}
        The corresponding simple roots are:
        $$
        \Delta^{(N)}=\{\epsilon_{1}+\epsilon_2\}\cup\{-\epsilon_i+\epsilon_{i+1}|\,1\leq i<N\},
        $$
        and the half sum of positive roots is:
        $$
        \rho^{(N)}=\frac{1}{2}\sum_{\alpha\in \Phi^{(N)}_{+}}\alpha=\sum_{1< i\leq N}(i-1)\epsilon_i.
        $$
    }
    \par{
        Now we show that $M_r$ is a generalized Verma module for the Lie algebra $\mathfrak{so}_{\infty}$. Define:
        \begin{align*}
        &\mathfrak{n}^{(N)}_{-}\stackrel{def.}{=}{\rm span}\{e_ke_l|\,k,l<0\},\,\mathfrak{l}^{(N)}\stackrel{def.}{=}{\rm span}\{e_ke_l|\,k<0,l>0\},\\
        &\mathfrak{u}^{(N)}\stackrel{def.}{=}{\rm span}\{e_ke_l|\,k,l>0\},\,\mathfrak{p}^{(N)}\stackrel{def.}{=}\mathfrak{l}^{(N)}\oplus \mathfrak{u}^{(N)}.
        \end{align*}
        Then
        \begin{align*}
        &\mathfrak{g}^{(N)}=\mathfrak{p}^{(N)}\oplus\mathfrak{n}^{(N)}_{-}=\mathfrak{l}^{(N)}\oplus\mathfrak{u}^{(N)}\oplus\mathfrak{n}^{(N)}_{-}.
        \end{align*}
        We also define $\Phi^{(N)}_{I}$ as follows
        $$
        \Phi^{(N)}_{I}\stackrel{def.}{=}\{-\epsilon_i+\epsilon_j|\,i<j\}\cup\{\epsilon_i-\epsilon_j|\,i<j\},
        $$
        then $\mathfrak{l}^{(N)}$ is a direct sum of $\mathfrak{h}^{(N)}$ and root spaces $(\mathfrak{g}^{(N)})_{\alpha}$, where $\alpha\in \Phi^{(N)}_{I}$.
    }
    \par{
        Define the 1-dimensional $\mathfrak{p}^{(N)}$-module of weight $\lambda^{(N)}\in (\mathfrak{h}^{(N)})^{*}$ spanned by the element $1$ such that:
        \begin{align*}
        &x\cdot 1=0,\quad h\cdot 1=\lambda^{(N)}(h)\cdot 1
         \quad \forall h\in \mathfrak{h}^{(N)}
         ,x\in (\mathfrak{g}^{(N)})_{\alpha},\alpha\in\Phi^{(N)}_{I}\cup \Phi^{(N)}_{+}.
        \end{align*}
        Then we obtain the following generalized Verma module $M_{I}(\lambda^{(N)})$:
        \begin{align}
        M_{I}(\lambda^{(N)})\stackrel{def.}{=}U(\mathfrak{g}^{(N)})\otimes_{U(\mathfrak{p}^{(N)})}\mathbb{C}\cdot 1\simeq U(\mathfrak{n}^{(N)}_{-})\cdot 1.\label{Mil}
        \end{align}
        The module $M_{I}(\lambda^{(N)})$ is a `generalized Verma module of scalar type' (see for example \cite{Hum}), and by (\ref{action2}), it is easy to calculate $\lambda^{(N)}$ explicitly:
        $$
            \lambda^{(N)}=-\frac{r}{2}\sum_{i=1,\cdots, N}\epsilon_i.
        $$
        We note that have an embedding of $\mathfrak{g}^{(N)}$-modules
        $$
            M_I({\lambda^{(N)}})\hookrightarrow M_r,
        $$
        and we have an exhaustive filtration:
        \begin{align}
        \{0\}\subseteq M_I(\lambda^{{(1)}})\subseteq \cdots \subseteq M_I({\lambda^{(N)}})\subseteq \cdots \subseteq M_r. \label{Mfilt}
        \end{align}
        Hence $M_r=V_{\mathcal{J},r}$ is a `generalized Verma module of scalar type' for $\mathfrak{so}_{\infty}$.
    }
    \par{
    We recall the following lemma on the simplicity of the scalar type generalized Verma module. This lemma is essentially due to Jantzen etc. \ref{}, and the reader can consult \cite{Hum}, Theorem 9.12 (a):
    \begin{lemma}
        If $\lambda^{(N)}$ is a weight for $\mathfrak{g}^{(N)}$ and
        \begin{align}
            \langle\lambda^{(N)}+\rho^{(N)},\beta^{\vee}\rangle\notin \mathbb{Z}_{>0},\quad\forall \beta\in\Phi^{(N)}_{+}-\Phi^{(N)}_{I}, \label{cond}
        \end{align}
        then $M_{I}(\lambda^{(N)})$ is a simple $\mathfrak{g}^{(N)}$-module. Conversely, if $M_{I}(\lambda^{(N)})$ is simple and $\lambda^{(N)}$ is regular, that is
        $$
            \langle\lambda^{(N)} +\rho^{(N)},\beta^{\vee}\rangle\neq 0
        $$
        for all $\beta\in \Delta$, then (\ref{cond}) also holds.
    \end{lemma}
    It is easy to compute that
    $$
        \beta^{\vee}=
        e_{-k}e_k+e_{-l}e_l
    $$
    for $\beta=\epsilon_k+\epsilon_l\in \Phi^{(N)}_{+}-\Phi^{(N)}_{I}=\{\epsilon_i+\epsilon_j|\,1\leq i<j\leq N\}$. Hence
    $$
        \langle \lambda^{(N)}+\rho^{(N)},\beta^{\vee}\rangle=
            -r+k+l.
    $$
    It is obvious that if $r\notin \mathbb{Z}$, then
    $$
        \langle \lambda^{(N)}+\rho^{(N)},\beta^{\vee}\rangle\notin\mathbb{Z}_{>0},\,\forall \beta\in\Phi^{(N)}_{+}-\Phi^{(N)}_{I}.
    $$
    By Lemma 2.4 $M_{I}(\lambda^{(N)})$ is simple as a $\mathfrak{g}^{(N)}$-module if $r\notin \mathbb{Z}$.
    }
    \par{
    We now conclude that $V_{\mathcal{J},r}$ is simple if $r\notin \mathbb{Z}$ by contradiction. Suppose the contrary that $V_{\mathcal{J},r}$ is not simple for $r\notin \mathbb{Z}$, then it has a proper $\mathfrak{so}_{\infty}$ submodule $M$. We deduce that $M\cap M_I(\lambda^{(N)})$ is a proper $\mathfrak{g}^{(N)}$-submodule of $M_I(\lambda^{(N)})$, for some $N$. This contradicts to the result that $M_I(\lambda^{(N)})$ is irreducible for all $N$ when $r\notin\mathbb{Z}$. Hence we conclude the proof in one direction.
    }
    \par{
        To prove the converse, we need to show that if $V_{\mathcal{J},r}$ is not simple, then $r\in\mathbb{Z}_{\neq 0}$. This happens if for some $N$, $M_I(\lambda^{(N)})$ is an reducible $\mathfrak{g}^{(N)}$-module. A direct calculation shows that for the simple root $\beta=\epsilon_1+\epsilon_2$ we have
        $$
            \langle \lambda^{(N)}+\rho^{(N)},\beta^{\vee}\rangle=
            -r+1,
        $$
        and
        $$
            \langle \lambda^{(N)}+\rho^{(N)},\beta^{\vee}\rangle=0
        $$
        for other simple roots $\beta$. Therefore $\lambda^{(N)}$ is regular if and only if $r\neq 1$. It is also easy to show that for $r\in \mathbb{Z}_{\neq 0,1}$, there exists an integer $N$ and $1\leq k,l\leq N$ such that
        $$
        \langle \lambda^{(N)}+\rho^{(N)},\beta^{\vee}\rangle=
            -r+k+l\in\mathbb{Z}_{>0},
        $$
        therefore by the second half of Lemma 4, $V_{\mathcal{J},r}$ is not simple when $r\in\mathbb{Z}_{\neq 0,1}$.
    }
    \par{
        It remains to show that $V_{\mathcal{J},1}$ is not simple. This follows from the fact that the VOA $\bar{V}_{\mathcal{J},1}$, is a quotient of $V_{\mathcal{J},1}$, and
        $$
            V_{\mathcal{J},1}\neq \bar{V}_{\mathcal{J},1}
        $$
        by comparing the dimensions of the corresponding graded vector spaces. Thus we conclude the proof of (3) in Theorem 1.
    }
    \par{
        It is interesting to know the explicit construction of the simple quotients of $V_{\mathcal{J},r}$ for all $r\in \mathbb{Z}_{\neq 0}$. This can be done using dual-pair type constructions similar to the one given in \cite{Z2}, and it will be given in a separate paper.
    }

    \par{
\textsc{School of Mathematics (Zhuhai), Sun Yat Sen University, Zhuhai, China.}\\
}
\par{
\textit{Email Address}: \textbf{hzhaoab@connect.ust.hk}
}
\end{document}